\documentclass[a4paper,conference]{IEEEtran}

\usepackage[pdftex]{graphicx}
\usepackage{array}
\usepackage{mdwmath}
\usepackage{mdwtab}
\usepackage{amssymb,latexsym}
\usepackage{stfloats}
\usepackage{amsmath}
\usepackage{subfig}
\usepackage{url}

\usepackage{algcompatible}

\usepackage{algorithm}

\usepackage[font={footnotesize}]{caption}

\usepackage{balance}

\usepackage[utf8]{inputenc}

\usepackage[T1]{fontenc}

\usepackage[compatible]{algpseudocode}

\renewcommand{\figurename}{Figure}

\renewcommand\IEEEkeywordsname{Keywords}

\makeatletter
\def\ps@IEEEtitlepagestyle{%
  \def\@oddfoot{\mycopyrightnotice}%
  \def\@evenfoot{}%
}
\def\mycopyrightnotice{
  {\footnotesize 979-8-3503-7542-8/24/\$31.00 ©2024 IEEE\hfill} 
  \gdef\mycopyrightnotice{}% just in case
}

\begin{document}

\title{Impact Analysis of Data Drift Towards The Development of Safety-Critical Automotive System}

\author{\IEEEauthorblockN{
Md Shahi Amran Hossain,
Abu Shad Ahammed,
Divya Prakash Biswas,
Roman Obermaisser}
\IEEEauthorblockA{University of Siegen
% University of Siegen\\
% Siegen, Germany}
% \IEEEauthorblockA{\IEEEauthorrefmark{2}
% University of Siegen\\
% Siegen, Germany}
% \IEEEauthorblockA{\IEEEauthorrefmark{3}
% University of Siegen\\
% Siegen, Germany}
% \IEEEauthorblockA{\IEEEauthorrefmark{4}
\\
Siegen, Germany}

{\it md.hossain@uni-siegen.de}}

\maketitle

\begin{abstract}
A significant part of contemporary research in autonomous vehicles is dedicated to the development of safety critical systems where state-of-the-art artificial intelligence (AI) algorithms, like computer vision (CV), can play a major role. Vision models have great potential for the real-time detection of numerous traffic signs and obstacles, which is essential to avoid accidents and protect human lives. Despite vast potential, computer vision-based systems have critical safety concerns too if the traffic condition drifts over time. This paper represents an analysis of how data drift can affect the performance of vision models in terms of traffic sign detection. The novelty in this research is provided through a YOLO-based fusion model that is trained with drifted data from the CARLA simulator and delivers a robust and enhanced performance in object detection. The enhanced model showed an average precision of 97.5\% compared to the 58.27\% precision of the original model. A detailed performance review of the original and fusion models is depicted in the paper, which promises to have a significant impact on safety-critical automotive systems. 

\end{abstract}

\begin{IEEEkeywords}
CARLA, Data drift, Automotive AI, Safety critical system, YOLO, Object detection
\end{IEEEkeywords}

\IEEEpeerreviewmaketitle

\section{Introduction}
\label{int}
The impact of artificial intelligence (AI) on self-driving vehicles, also known as autonomous vehicles, is profound and multifaceted, defining not only how vehicles should operate safely, but also providing an infrastructure to incorporate these vehicles into a broader transportation system. The presence of AI is growing day by day in safety-critical systems, particularly in the autonomous driving domain, to perform tasks such as object detection, lane maintenance, lane switching, and navigation towards a destination \cite{fan2021real}. Such technological advancements that are expected to revolutionize the future of autonomous vehicles are mostly driven by computer vision algorithms, a cornerstone of automotive AI that enables computers to interpret and understand visual information from the world. Already, industrial giants like NVIDIA, Mobileye, Intel, Qualcomm, and Samsung have made significant progress in developing real-time vision applications for autonomous vehicles \cite{lecun2015deep}. Although vision algorithms have high potential, they diverge in many aspects like lack of transparency, nondeterministic outcomes, reliance on a subset of input domain, instability, etc., and all these characteristics are diametric to the safety-critical system. In automotive sector, the major shortcomings are it's limited robustness i.e. poor object detection and recognition capabilities in diverse conditions, which can hinder safe driving at road.\\
However, it is necessary to use vision-based AI algorithms in the autonomous driving domain to meet current demands like improved road safety, traffic efficiency, and reduction in human errors leading to fewer accidents \cite{azam2020system}. The study depicted in this research here highlights the limitation of the vision model, specifically its reduced performance in detecting objects such as road signs under challenging traffic conditions. There are several reasons why a computer vision model may have distorted performance, and one of them is lack of proper training. If there is insufficient, imbalanced or noisy information in the data set, or the data set contains low diversity, which results in high data drift in real world simulation, then the performance of vision model will not be as expected. Data drift occurs when new types of scenario appear in the environment that the model has not been trained to recognize.\\
The study proposes a solution through the development of two object detection models utilizing YOLO version 8, renowned for its precision and minimal inference time. The models have been constructed based on traffic data that encompass eight different road signs: Round-shaped 30, 60 and 90 kilometers speed limit;  Square shaped 30, 60 and 90 kilometers speed limit; Stop sign. The first model is developed with road sign images only in normal environmental and traffic conditions. The second model is the fusion model, which was developed with road sign images in both normal and drifted environments. The aim was to verify how the standard model and the fusion model perform when an identical set of drifted datasets is provided for test purposes to both models. The image data was generated using various road patterns within CARLA, an automotive simulation environment. The paper also includes a comparative performance analysis to demonstrate the improved performance of the fusion model in terms of selected evaluation criteria.\\
The remainder of the manuscript is organized as follows. Section \ref{lit} provides a meta-analysis of the existing research works, methods and techniques by different authors. Section \ref{ai} will discuss the prospect of using AI in automotive safety critical system. Section \ref{CARLA} will discuss briefly the CARLA simulator and how it contributes to the automotive industry. A conceptual overview of data drift with possible examples will be shown in section \ref{drift}. Section \ref{org} discusses the organization and management of the normal and drifted data with preprocessing steps. The development of the YOLO model and selected hyperparameters is discussed in Section \ref{obj}. Section \ref{res} shows the complete performance analysis of the two object detection models to reflect the validity of this investigation. Finally, Section \ref{conc} provides a conclusion, mentioning the existing challenges and a look at the possibilities in the future.
\section{Literature Review}  
\label{lit}
The research we conducted is focused on analyzing the impact of data drift on the critical computer vision-based automotive safety system. We also provided our research solution with a fusion model that used drifted data generated from the CARLA simulator to improve object detection performance. Even though prior studies have examined the effects of data drift in vision models, our research uniquely investigates its implications in real-time automotive applications.\\
Artificial intelligence is a fundamental component currently employed in designing autonomous cars as it can act like the brain of self-driving cars, enabling them to carry out tasks such as real-time object detection, lane keeping or switching, and navigation towards a destination following safety regulations. As referred by Bathla et al. \cite{bathla2022autonomous}, AI technology can completely replace humans with intelligent automation while offering better safety and intelligent movement of vehicles. But there are also formidable challenges in regards to safety and reliability when employing AI algorithms to autonomous vehicles across diverse driving scenarios \cite{fagnant2015preparing}. As suggested by \cite{kannally2023human}, the integration of AI in automotive systems necessitates robust human-AI teaming to ensure effective collaboration and decision-making between AI technologies and human operators. One of the main reasons identified as a risk is data drift over time or in a changed scenario which degrade the performance of AI in real time. Several important studies \cite{vela2022temporal} \cite{ackerman2021automatically} \cite{ackerman2020detection} have been identified that explain how data drift in automotive AI can greatly affect the performance and dependability of AI models. Data drift is explained as the situation where the distribution of new data evolves over time, causing a discrepancy between the data used for training and the data encountered during operation, which can result in a decline in the model's performance and reliability. As our research is focused on change impact analysis of data drift, we used a synthetic data set generated from CARLA simulator that can replicate real-world scenarios accurately as described in \cite{cao2023data} \cite{rosende2023urban}.\\
\section{AI in Automotive Safety Critical Systems}
\label{ai}
An automotive safety critical system encompasses any component, subsystem, or system within a vehicle that aims to safeguard the safety of passengers, pedestrians, and other road users. Modern automotive systems consist of different types of microprocessors and microcontrollers that vary in functionality, as well as their importance in ensuring safety critical systems. Recent developments in autonomous cars have increased the need for better safety standards, and a lot of research is still ongoing due to the increasing complexity and interconnectivity of modern cars. These systems e.g. anti-lock breaking system (ABS), airbag, adaptive cruise control, etc. are designed to perform their function or set of functions under all circumstances, including when certain parts of the vehicle fail or when the vehicle is operating under extreme conditions. Artificial Intelligence (AI), as one of the most talked topics now a days, can also play a vital role in designing or updating such systems. It can continuously monitor automotive sensor data for anomalies that may indicate potential safety risks, such as sudden changes in road users' behavior that may lead to crush,distorted road signs, or bizarre weather conditions increasing difficulty to detect the road signs correctly. By predicting anomalies in real time, these systems can alert drivers or take corrective actions autonomously or semi-autonomously to prevent accidents.\\
The computer vision (CV) algorithm employed in our research is essential for the development of autonomous vehicles, as it supports vehicles to detect and recognize various elements of their driving path, including static objects like traffic lights, signs and crossings, as well as dynamic objects like other vehicles, pedestrians, and cyclists \cite{Zepeda2022} \cite{mahima2021adversarial}. Although vision algorithms have great potential in the automotive sector, as per \cite{amodei2016concrete}, we need to seriously consider the societal impacts of CV: the problem of preventing accidents using vision or AI algorithms. Because CV models are not uncommon to have poor performance due to changes in lighting conditions, adverse weather, and occlusions \cite{rezaei2015robust}. Another issue is the sensitivity of CV algorithms to image quality as video images extracted in real time from practical scenarios can be distorted and drifted, which poses a challenge in maintaining consistent performance. Even slight variations in image quality can result in significant declines in algorithm effectiveness, emphasizing the need for robust algorithms capable of adjusting to different environmental and traffic conditions \cite{yahiaoui2019overview}.
\section{CARLA Simulator}
\label{CARLA}
Car Learning to ACT (CARLA) is an open-source autonomous driving simulator born from a joint initiative of researchers from Intel Labs, Toyota Research Institute, and the Computer Vision Center of Barcelona and was released in 2017 \cite{ly2020learning}. The simulator has been developed from the ground up to support the development, training, and validation of autonomous driving systems. In particular, it provides a flexible and highly realistic simulation environment that features a variety of weather conditions, different types of automotive sensors, accurately modeled vehicles, and a diverse set of environments. The major functionalities of CARLA simulators are listed below.
\begin{itemize}
    \item Open-source, flexible, and scalable platform that allows multiple clients to manage various actors within the simulation environment. Users can effortlessly control any aspect related to simulation, such as traffic generation, road signs manipulation, pedestrian behavior, weather conditions, sensors, and more, using a Python API.
    \item Offers rich and realistic simulation environments with high fidelity graphics using Unreal Engine that provides developers access to a wide range of urban settings, such as cities, rural roads, and highways, with various traffic situations, pedestrians, and vehicles.
    \item Wide range of automotive sensors such as Camera, Lidar, Radar, GPS, and IMU are available in CARLA to configure and simulate based on user preference
    \item Predefined and custom scenarios can be created to systematically test and validate driving policies under controlled but realistic conditions. CARLA provides tools like Scenario Definition Language (SDL), ASAM OpenScenario, or RoadRunner to create and test such scenarios
    \item Predefined and custom scenarios can be created to systematically test and validate driving policies under controlled but realistic conditions
    \item CARLA can even be programmed for the modeling of AI agents with varying levels of autonomy to interact with and adapt to the simulated world
\end{itemize}
The impact of CARLA in automotive industry is very high in terms of cost-effective test simulation and scenario generation \cite{won2022verification}. Detecting objects like cones in construction zone scenarios, road signs or pedestrian identification for autonomous driving is a core expectation that can be tested using CARLA. In our research, we have used CARLA to create synthetic data set that represents drifts in normal traffic signs like speed limiters or STOP. Custom mapping option of CARLA was also utilized to modify objects and create drift scenarios. 
\section{Data Drift}
\label{drift}
\begin{figure}[!b]
\center
\subfloat[Normal Road Sign]{\label{fig:std1}
\includegraphics[width=0.24\textwidth]{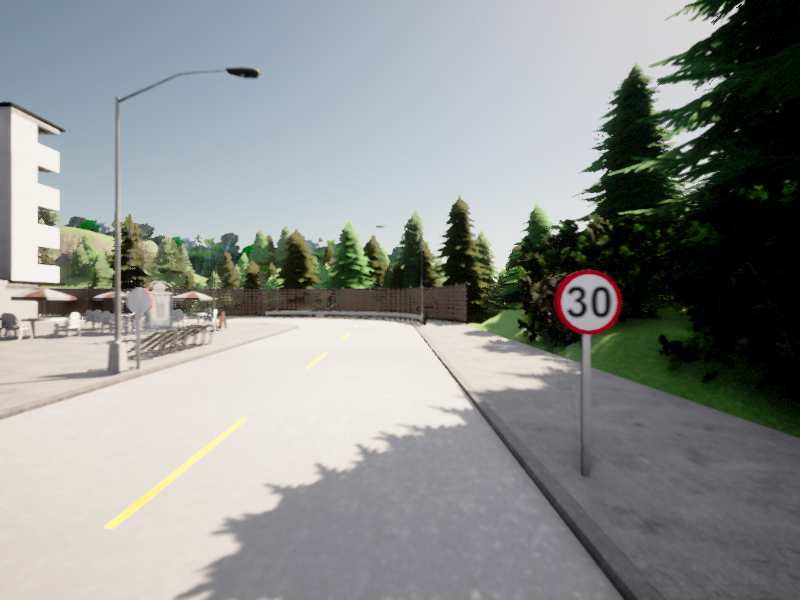}}
\subfloat[Tilted Road Sign]{\label{fig:drift1}
\includegraphics[width=0.24\textwidth]{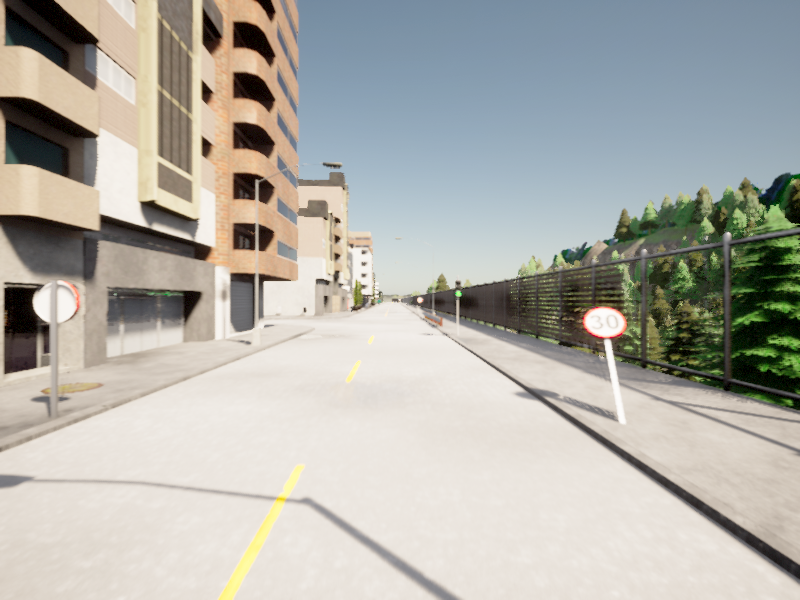}}
\\
\subfloat[Normal Sunlight]{\label{fig:std2}
\includegraphics[width=0.24\textwidth]{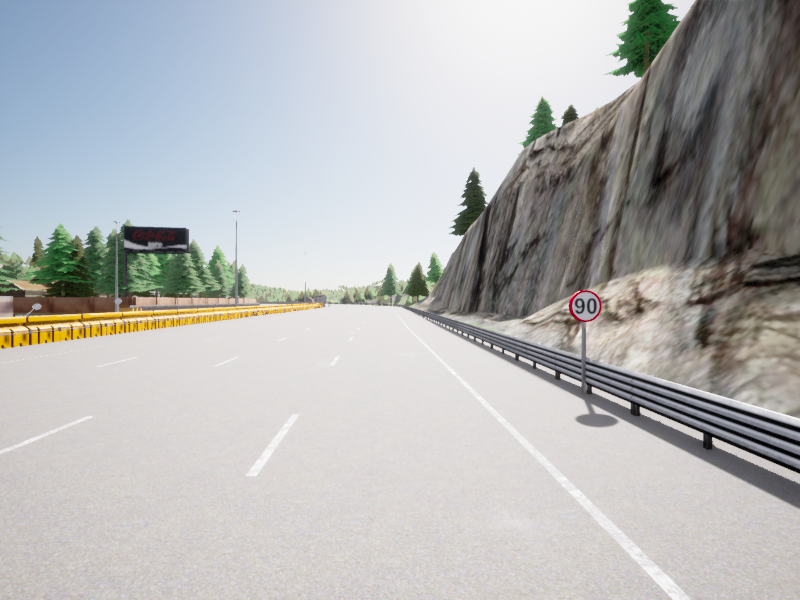}}
\subfloat[Low Light at Night]{\label{fig:drift2}
\includegraphics[width=0.24\textwidth]{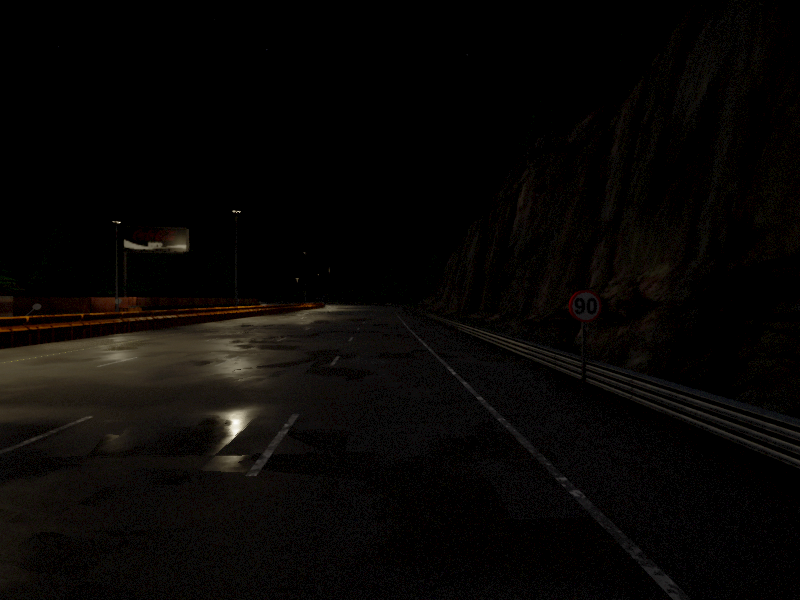}}
\caption{ An Illustration of Data Drift under Normal and Altered Environmental Conditions}
\label{fig:drift}
\end{figure}

Data drift is one of the major challenges for deploying computer vision models in a fluid environment. It implies changes in the statistical properties and characteristics over time of the input domain which may have a significant effect on the performance of the machine learning model. In safety-critical applications of autonomous vehicles, aerospace, or nuclear power plants where reliability and robustness are the main priority, data drift is a primary concern. Traditionally, the training data set that is used to train an CV model is historical data from reliable sources.  The CV model is expected to perform well if the properties of the test input data are close to those of that training data set. But it can also happen that the trained model when deployed in the real world gets data from the practical scenario which is entirely different due to the adapting nature of the environment. In that case, the train data set may become antiquated, and the model's accuracy might drop dramatically. As an example, if a traffic sign is tilted more than its natural position which is not considered during training period may reflects more sunlight as shown in Figure \ref{fig:std1} and \ref{fig:drift1}. Similarly, there can be conditions like poor night light which may make it difficult for autonomous cars to detect the road signs as depicted in Figure \ref{fig:std2} and \ref{fig:drift2}.
In addition, for the highlighted drifts shown in Figure \ref{fig:drift}, there are other types of data drift available in the domain of the traffic sign recognition models. The drifts considered in our research are illumination drift, weather drift, seasonal drift, or camera drift. As explained in Section \ref{CARLA}, the CARLA simulator with its flexible map customization function was used to create all these drift-based custom scenarios. A total of 2,192 drifted image sets representing all above-mentioned drift scenarios have been generated that can be used to train the fusion model.
\section{Data Organization and Pre-processing}
\label{org}
The organization and pre-processing of the date set is a critical step in the computer vision pipeline as it directly impacts the performance, accuracy, and efficiency of the object detection models.  The data set planned in our research consist of standard and drifted images of traffic scenarios later to be used for YOLO models. It comprises eight traffic signs as described in Section \ref{int}, with all images sourced from the CARLA simulator.  It was also possible to collect these image data with a camera from real road environment. But this process is time-consuming, requires a lot of manpower, and on top involves several risk factors compared to our way of utilizing the CARLA simulator, as it is an open source environment and provides more flexibility to support the development of autonomous driving systems with digital assets e.g. different traffic signs, urban and rural road layouts, different types of dynamic and static actors, etc. Initially, we have created two data sets: one comprises solely 2,017 images without any data drift, while the other includes 5,063 images consist of both drifted and normal images. Data augmentation techniques have been applied to increase the data volume and make it balance. In the process of preparing this data set, a few steps are taken as follows:
\begin{figure}[!b]
\center
\includegraphics[width=0.28\textwidth]{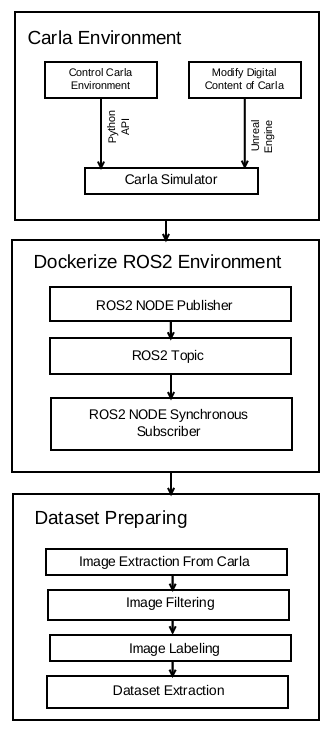}
\caption{An Overview of Data Organization and Pre-processing}
\label{fig:sample} 
\end{figure}
\begin{itemize}
\item \textbf{Data Collection}: This is the first step where the CARLA simulator environment has been properly set up. For standard images, different maps have been chosen for getting the diverse environment, and the weather conditions such as the sun position, fog density wetness, etc., have been kept in normal visibility range.  For drifted images, additional steps have been taken in comparison to the standard image collection procedure. First, the Sun position of the CARLA environment has been controlled through the Python API to get images of different times of day and night, which has a significant effect on images due to sunlight reflection and simulating the lack of sufficient daylight. Likewise, adverse weather conditions like rain or foggy environments have also been simulated through the Python API to ensure that the road sign images are highly drifted. The traffic signs that have been placed on all the maps are rotated at different angles through Unreal Engine 4 manually because these static actors cannot be controlled by the Python API. For collecting the images from the CARLA simulator, Robot Operating System 2 bridge have been used to maintain the standard communication. In the ROS2 bridge, there is a publisher node that is responsible for publishing image data to a ROS2 node continuously while a vehicle with an RGB camera is under driving conditions in the simulator environment. Afterwards, a synchronous subscriber node has been implemented to extract the images from that node as soon as it is saved. This whole process has been contained and maintained through Docker to avoid any dependency conflicts.
\item \textbf{Data Selection}: After completing the collection process, images were chosen manually which included desired traffic signs, and all others without signs were discarded. The next step was to check the quality of the collected images which is necessary to create a quality model training pipeline.
\item \textbf{Data Annotation and Augmentation}: In this step, the collected images have been annotated using computer vision annotation tool (CVAT) with their respective classes in YOLO format. CVAT, developed by Intel, is a free, open source, web-based image and video annotation tool used for labeling data for computer vision algorithms. Subsequently, techniques for data augmentation, including mirroring, blurring, and rotating, were employed to ensure the dataset is balanced and exhibits adequate diversity. 
\item \textbf{Data Splitting}: The last step in data organization was to split the images into three sets: training, testing, and validation. For the original model, 1,613 images were considered for training and 404 images for validation. The fusion model had 4,296 images for training and 767 images for validation. The training images considered have adequate variations to train the model for all environment aspects. For testing purpose of both models, four identical samples of each of the 8 road signs were considered with variation in traffic condition.
\end{itemize}
\section{Object Detection Model}
\label{obj}
You Only Look Once (YOLO) is an open source family of deep learning-based object detection models, designed to detect and classify objects in real-time from images and videos. There are a few object detection algorithms available such as R-CNN, SSD (Single Shot Multibox Detector), and Retinate. However, when it comes to speed and efficiency, YOLO surpasses all other algorithms, particularly its most recent versions. In the automotive sector, where autonomous driving demands rapid decision-making for navigation, the YOLO model is an ideal fit. Moreover, YOLO tends to generalize well to new images and environments because it considers the entire image and learns more generalized features during training compared to the other object detection algorithms.\\
For this study, YOLO version 8 has been chosen, which is built on cutting-edge advancements in deep learning and computer vision and offers the best performance in speed and accuracy \cite{YOLO}. This version of YOLO also typically outperformed its next generation (version 9) in both precision and sensitivity \cite{Encord2024}. As the first step in developing the YOLO based object detection model, the training phase was initiated with the standard data set generated from the CARLA simulator environment. According to the Ultralytics YOLO documentation \cite{UltralyticsTrain2024}, the split images of the train test validation were saved in their respective folders with a 'yaml' file that contains information regarding the annotated classes and the path to each data directory. Python version 3.11 was used to create the environment necessary to run the codes. For this training, the following hyperparameters were used:
\begin{itemize}
    \item Optimizer: Adamw
    \item Learning rate: 0.000833
    \item Epochs: 100
    \item Confidence score: 0.2
\end{itemize}
\begin{figure}[!h]
\center
\subfloat[Training with Standard  Dataset]{\label{fig:first}
\includegraphics[width=0.20\textwidth]{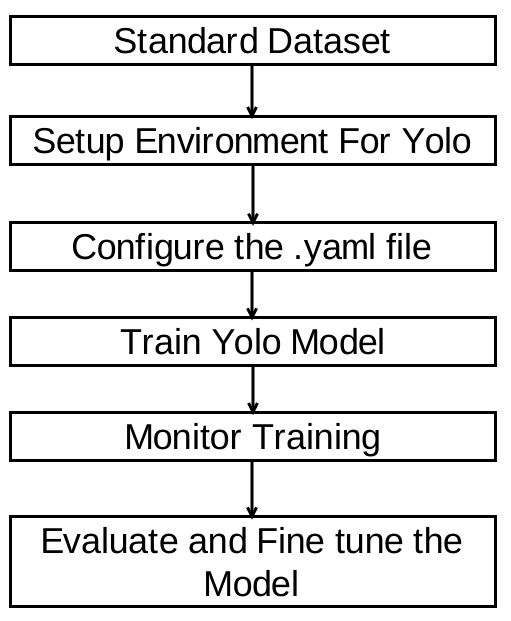}}
\subfloat[Training with Drifted  Dataset]{\label{fig:second}
\includegraphics[width=0.28\textwidth]{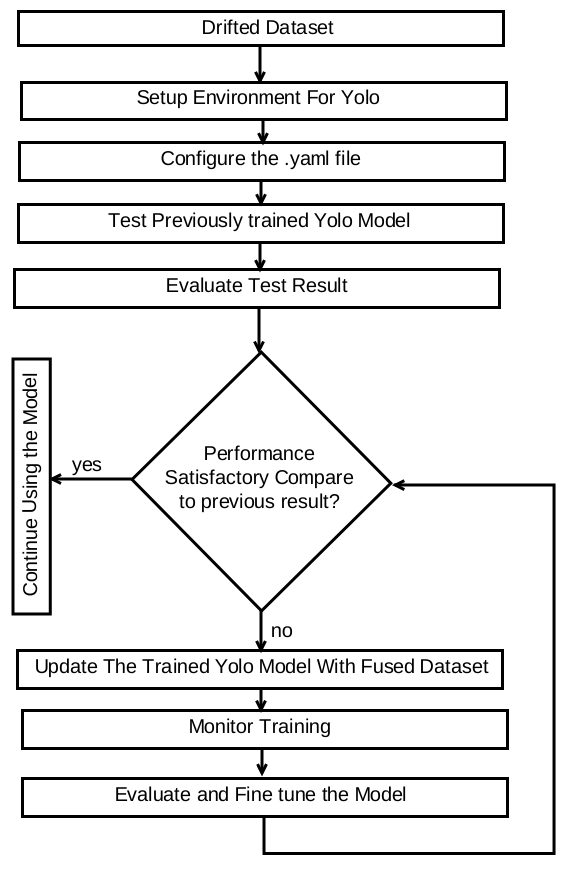}}
\caption{YOLO training Process for Standard and Drifted Dataset.}
\label{fig:trg}
\end{figure}
A complete overview of the training process can be seen from Figure \ref{fig:trg}.
The training steps were executed two times, one with the standard data set and then with the combination of drifted and standard datasets. The ultimate goal was to apply the test data set with the trained models to determine how data drift affects the models' predictive performance.
\section{Results and Discussion} 
\label{res}
\begin{figure}[!b]
\center
\includegraphics[width=0.52\textwidth]{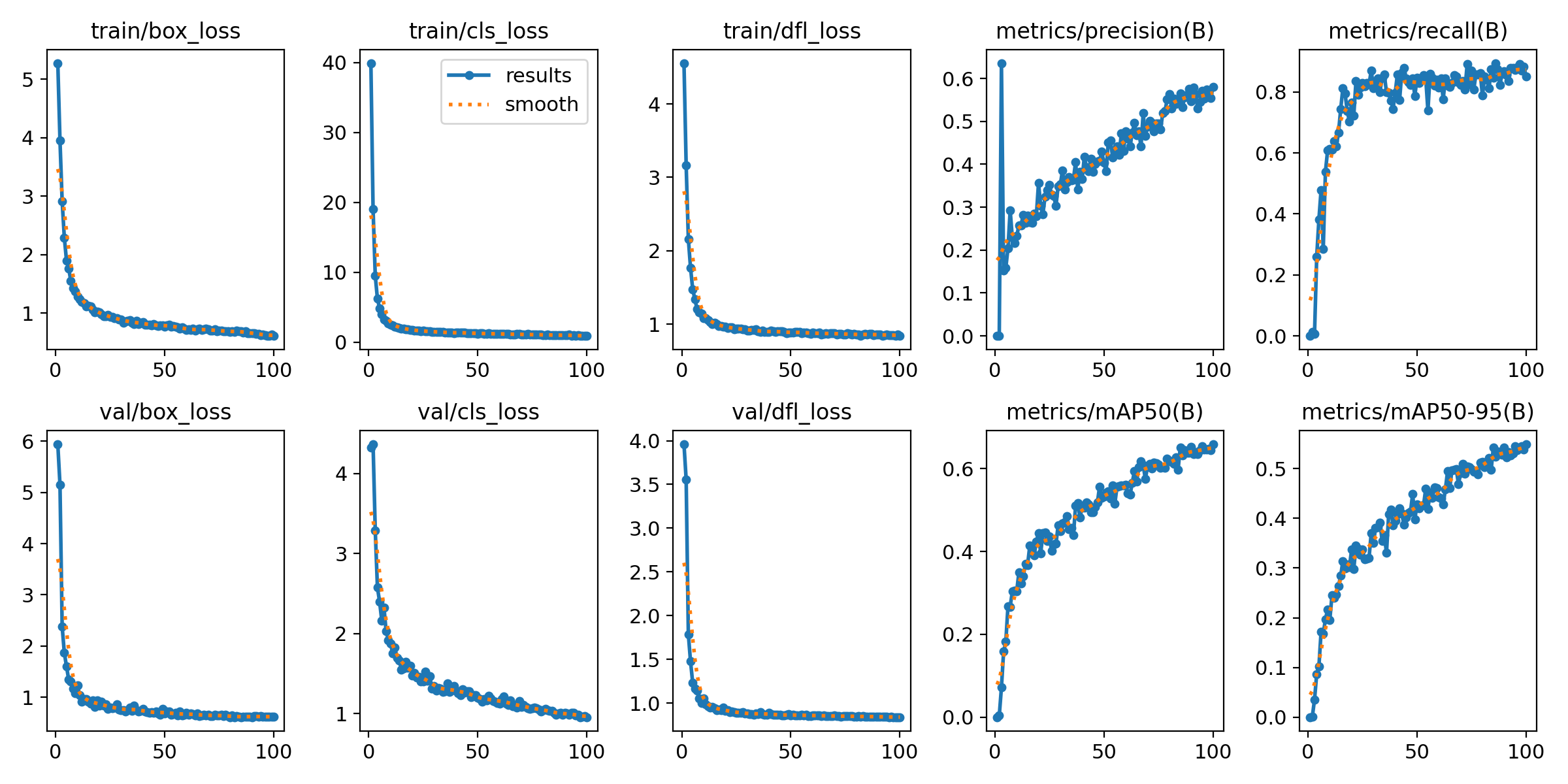}
\caption{Results Generated from Original Model}
\label{fig:OrgRes} 
\end{figure}

\begin{figure}[!b]
\center
\includegraphics[width=0.52\textwidth]{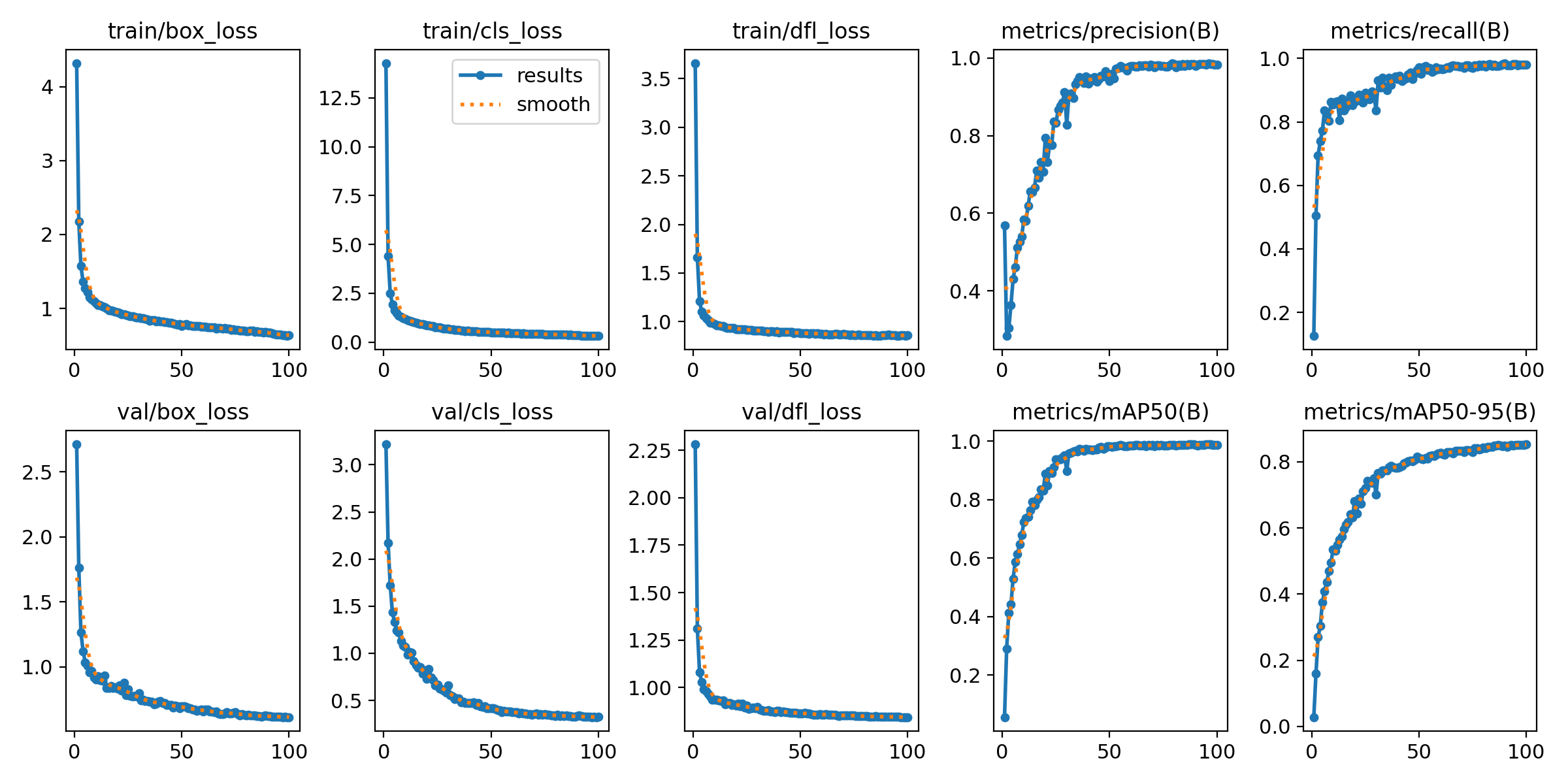}
\caption{Results Generated from Fusion Model}
\label{fig:FusRes} 
\end{figure}

Some basic metrics were used to assess the performance of the YOLO model in recognizing and classifying traffic signs as below \cite{quach2023evaluating}:
\begin{itemize}
\item Precision: Measures the accuracy of positive prediction of a model. Indicates the model’s ability to avoid false alarms
\item Recall: Represents the ability to identify all relevant positive instances in the data set
\item F1 Score: Provides a balance between the measure of precision and recall which gives insight into both the correctness of the positive prediction and the ability to find all the positive instances
\item Mean Average Precision 50 (mAP50): Evaluates how proficiently the model can detect objects. It calculates the average precision-recall curves across all classes at a 50\% Intersection over Union (IoU) threshold.
\item Mean Average Precision 50-95 (mAP50-95): Provides a more thorough and demanding measure of the model’s ability. Instead of 50\%, the IoU threshold is set in the range of 50-95\% to calculate the average precision recall curves in all classes.
\end{itemize}
After each training stage, the YOLO model typically generates some performance curves like confusion metrics, precision recall, F1 score, training and validation loss, etc. that provides insights into various aspects of the model's accuracy, efficiency, and generalization ability. In Figures \ref{fig:OrgRes} and \ref{fig:FusRes}, autogenerated curves are shown to understand the behavior and performance of the model.
As illustrated in the figures, the models demonstrate symmetrical advancement in both loss curves and performance metrics throughout each training epoch. The validation loss curves for both models showed a decreasing trend, which indicates that the model is not overfitting and generalizing well to the validation data. The models also showed high precision and recall values that represent a low false positive and negative rate. The mAP50 curves suggest that the models are effectively identifying and pinpointing objects with a 50\% intersection over union with the actual bounding boxes. In addition, the map50-95 curves showed that the models performed well across various levels of localization precision.\\
The next step was to measure the performance of the models with the test dataset that contains only drifted data. To achieve this, we selected the top performing model based on epoch-wise performance, which was automatically generated during training. The test sequence was run in prediction mode using the same hyperparameters as the training phase. A summary of the performance of the models is shown in Tables \ref{tab:org} and \ref{tab:fused}.
\begin{table}[!h]
\renewcommand{\arraystretch}{1.3}
\caption{\textsc{Validation and Test Result of Original Model}}
\label{tab:org}
\centering
\begin{tabular}{|l|c|c|}
\hline
\textbf{Criteria} & \textbf{Validation} & \textbf{Test} \\ \hline
Precision & 0.9899 & 0.5827 \\ \hline
Recall & 0.9920 & 0.5562  \\ \hline
F1-Score & 0.9909 & 0.5691  \\ \hline
mAP50 & 0.9891 & 0.3798  \\ \hline
mAP50-95 & 0.9126 & 0.4000  \\ \hline
\end{tabular}
\end{table}
\begin{table}[!h]
\renewcommand{\arraystretch}{1.3}
\caption{\textsc{Validation and Test Result of Fusion Model}}
\label{tab:fused}
\centering
\begin{tabular}{|l|c|c|}
\hline
\textbf{Criteria} & \textbf{Validation} & \textbf{Test} \\ \hline
Precision & 0.9827 & 0.9750 \\ \hline
Recall & 0.9799 & 0.9687  \\ \hline
F1-Score & 0.9812 & 0.9718  \\ \hline
mAP50 & 0.9872 & 0.9800  \\ \hline
mAP50-95 & 0.8700 & 0.7716  \\ \hline
\end{tabular}
\end{table}
The tables illustrate that while the models' validation performance during training was nearly identical, there was a significant difference when the test data set was applied. The precision, recall and F1 score values for the original model are near 50\% only and the map values are well below standard. It indicates that due to the introduction of data drift in the test images, the original model was unable to keep up with its expected performance. However, the fusion model showed much better performance in terms of test data. There was a small drop in performance for the validation set, but it is acceptable. The test results demonstrated a score of more than 95\% for precision, recall and F1 which is approximately 40\% better than the original model. This validates our approach of enhancing the accuracy of the object detection model by integrating drifted data from the CARLA simulator.
\section{Conclusion}
\label{conc}
Our research provided a thorough analysis on how data drift may impact the performance of traffic sign detection model developed with computer vision algorithm. The adverse effects of data drift were mitigated by training the detection model with a drifted dataset generated by CARLA, leading to enhanced safety critical automotive systems. This also highlights the future potential of AI models usage while complying safety standards like ISO26262 and ISO8800 in automotive industry.
\section*{Acknowledgment}
This work has been funded by the Federal Ministry of Education and Research (BMBF) as
part of AutoDevSafeOps (01IS22087Q).
\bibliographystyle{ieeetr}
\bibliography{reference}

\end{document}